\documentclass{article}

\usepackage[OT2,T1]{fontenc}
\DeclareSymbolFont{cyrletters}{OT2}{wncyr}{m}{n}
\DeclareMathSymbol{\Sha}{\mathalpha}{cyrletters}{"58}

\usepackage{amsbsy,amssymb,amscd,amsfonts,latexsym,amstext,delarray,
amsmath, diagrams}  \headsep=15pt
\usepackage[all,knot]{xy}
\def\piet{\pi_{1,\mbox{\'et}}}
\def\cP{ {\cal P} }
\def\Un{ \mbox{Un} }

\def\cX{{\cal X} }

\def\bX{ \bar{X}}
\def\Ext{{ \mbox{Ext} }}

\def\dim{{ \mbox{dim} }}
\def\Spec{{ \mbox{Spec} }}

\def\Hom{{ \mbox{Hom} }}

\def\ra{{ \rightarrow }}

\def\a{{ \alpha }}
\def\b{{ \beta} }

\def\d{{ \delta }}

\def\F{ {\mathbb F} }

\def\Un{ \mbox{Un} }
\def\Vect{ \mbox{Vect} }

\def\C{{ \mathbb{C} }}

\def\bs{ \backslash}

\def\G{{ \Gamma }}
\def\Gal{{ \mbox{Gal} }}

\def\bQ{\bar{\Q}}

\def\cC{ {\cal C}}
\def\Z{{ \mathbb{Z}}}
\def\cC{ {\cal C} }

\def\bq{\begin{quote}}
\def\eq{\end{quote}}
\def\Aut{ \mbox{Aut}}
\def\Isom{ \mbox{Isom} }

\newtheorem{thm}{Theorem}
\newtheorem{cor}[thm]{Corollary}

\def\Q{\mathbb{Q}}

\def\cV{ {\cal V} }

\def\P{ {\bf P}}

\def\Sym{\mbox{Sym}}
\def\be{\begin{equation}}
\def\ee{\end{equation}}
\def\k{ \kappa}

\def\cL{{\cal L}}

\def\cL{ {\cal L}}

\def\bE{ \bar{E}}

\def\loc{\mbox{loc}}

\def\Pet{\pi_1^{et}(\bX;b,x)}
\def\piet{\pi_1^{et}(\bX,b)}
\def\b{ \bar}

\def\L{\Lambda}

\def\cP{{\cal P}}

\def\a{\alpha}

\def\bd{\begin{diagram} }
\def\ed{ \end{diagram}}

\def\b{\beta}
\def\tb{\tilde{b}}

\def\bs{\begin{slide}}
\def\es{\end{slide} }

\def\Cov{\mbox{Cov}}

\def\bd{ \begin{diagram}}
\def\ed{\end{diagram} }

\def\Cov{\mbox{Cov}}

\def\Pux{\pi_1^{u,\Q_p}(\bX;b,x)}

\title{
Galois Theory and Diophantine geometry}
\author{Minhyong Kim}
\begin{document}
\maketitle
\begin{flushleft}
{\em Lecture at Cambridge workshop, July, 2009}
\end{flushleft}
\medskip

The author must confess to having contemplated for some years a diagram of the following sort.
$$\bd
& & \fbox{Diophantine geometry} & & \\
&\ruTo & & \luTo & \\
\fbox{theory of motives}& & & &\fbox{anabelian geometry}
\ed$$
To a large extent, the investigations to be brought up today arise from a curious
inadequacy having to do with the arrow on the left. On the one hand,
it is widely acknowledged that the theory of motives finds a strong source of
inspiration in Diophantine geometry, inasmuch so many of the
structures, conjectures, and results therein have as model the
conjecture of Birch and Swinnerton-Dyer, where the concern is
with rational points on elliptic curves that can be as simple as
$$x^3+y^3=1729.$$
Even in the general form discovered by Deligne, Beilinson, Bloch and Kato, (see, for example, \cite{motives1})
it is clear that motivic $L$-functions are supposed, in an ideal world, to
give access to invariants in arithmetic geometry of a {\em Diophantine nature}.
The difficulty arises when we focus on the very primitive concerns of Diophantine
geometry, which might broadly be characterized as the study of
maps between schemes of finite type over $\Z$ or $\Q$.
One might attempt, for example, to define the points of a motive $M$ over
$\Q$ using a formula like
$$\Ext^1(\Q(0),M)$$
or even
$$\mbox{RHom}(\Q(0), M),$$
hoping it eventually to be adequate in a large number of situations.
However, even in the best of all worlds, this formula will never provide
direct access to the points of a scheme, except in very special situations
like $M=H_1(A)$ with $A$ an abelian variety. This is a critical limitation
of the abelian nature  of motives, rendering it quite
difficult to find direct applications to any mildly non-abelian
Diophantine problem, say that posed by a curve of genus 2.
It is worth remarking that this limitation is essentially by design,
since the whole point of the motivic category is to {\em linearize}
by increasing the number of morphisms\footnote{Even then, we complain that
there are not enough.}. Of course we should pause to acknowledge the
role of technology that is more or less motivic in nature within two of the most celebrated
Diophantine results of our times, namely the theorems of Faltings
and of Wiles. But there, the idea is to constrain points on a non-abelian
variety by forcing  them to {\em parametrize} motives of a very special type. The method of achieving
this is highly ingenious in each case and, therefore, underscores our concern that
it is rather unlikely to be part of a general system, and certainly not
of the motivic philosophy as it stands.

Much has been written about the meaning of anabelian geometry, with a general tendency
to retreat to the realm of curves as the only firm ground on
which to venture real assertions or conjectures. We as well will  proceed to use $X$ to denote
a smooth projective curve of genus at least two over $\Q$. The basic anabelian
proposal then is to replace the Ext group that appeared above
by the topological space
$$H^1(G,\piet),$$
the non-abelian continuous cohomology \cite{serre1} of the absolute Galois group $G=\Gal(\bQ/\Q)$
of $\Q$ with coefficients in the profinite \'etale fundamental group
of $X$. The notation will suggest that a rational basepoint $b\in X(\Q)$ has
been introduced. Many anabelian results do not require it \cite{NTM}, but
the Diophantine issues discussed today will gain in clarity
by having it at the outset, even if the resulting restriction may
appear as serious to many.
An immediate relation to the full set of points is established by way of a non-abelian Albanese
map
$$ X(\Q) \rTo^{\k^{na}} H^1(G,\piet);$$
$$x\mapsto [\Pet].$$
 We remind ourselves that the
definition of fundamental groups in the style of Grothendieck \cite{szamuely}
typically starts from a suitable category over $X$, in this case that
of finite \'etale covers of $$\bX=X\times_{\Spec(\Q)}\Spec(\bQ)$$ that we might denote
by
$$\Cov(\bX).$$
The choice of any point $y\in \bX$ determines a fiber functor
$$F_y:\Cov(\bX) \rTo \mbox{Finite Sets},$$
using which the fundamental group is defined to be
$$\pi_1^{et}(\bX,y):=Aut(F_y),$$
in the sense of invertible natural transformations familiar from category
theory\footnote{The reader unfamiliar with
such notions would do well to think about  the case of a functor
$$F:\mathbb{N}^{op}\ra \cC$$
whose source is the category of natural numbers with a single morphism
from $n$ to $m$ for each pair $m\leq n$. Of course this is just a sequence
$$\ra F(3) \ra F(2) \ra F(1) \ra F(0)$$
of objects in $\cC$, and an automorphism of $F$ is a compatible sequence $(g_i)_{i\in \mathbb{N}}$
of automorphisms
$$g_i:F(i)\simeq F(i).$$ For a general
$F: {\cal B} \ra \cC$, it is profitable to think of ${\cal B}$ as a complicated
indexing set for things in ${\cal C}$.}. Given two points $y$ and $z$, there is also the set of
\'etale paths
$$\pi_1^{et}(\bX;y,z):=\Isom(F_y,F_z)$$
from $y$ to $z$ that the bare definitions equip with a right action of $\pi_1^{et}(\bX,y)$, turning it thereby into a torsor for the fundamental group. When $y$ and $z$ are rational points, the naturality of
the constructions equips all objects  with a compatible action of $G$, appearing in the
non-abelian cohomology set and the definition of the map
$\k^{na}$.

The context should make it clear that $H^1(G,\piet)$ can be understood as a non-abelian
Jacobian in an \'etale profinite realization, where the analogy might be strengthened by the
interpretation of the $G$-action as defining a sheaf on $\Spec(\Q)$ and
 $H^1(G,\piet)$   as the moduli space of
torsors for $\piet$ in the \'etale topos of $\Spec(\Q)$. It is instructive to compare
this space with the moduli space $Bun_n(X)$ of rank $n$
vectors bundles on $X$ for $n\geq 2$. Their study was initiated in a famous paper of Andr\'e Weil \cite{weil2}
whose title suggests the intention of the
author to regard them also as non-abelian Jacobians. Perhaps less well-known is the main motivation
of the paper, which the introduction essentially states to be the study of
rational points on curves of higher genus. Weil had at that point already expected
non-abelian fundamental groups to intervene somehow in a proof of the Mordell conjecture,
except that a reasonable arithmetic theory of $\pi_1$ was not available at the time. In order to make the connection to fields of definition, Weil
proceeded to
interpret the representations of the fundamental group in terms of algebraic
vector bundles, whose moduli would then have the same field of definition
as the curve. In this sense, the paper is very much a continuation of
Weil's thesis \cite{weil1}, where an algebraic interpretation of the Jacobian is attempted with the same
goal in mind, however with only the partial success noted by Hadamard. The spaces
$Bun_n(X)$ of course fared no better, and it is perhaps
sensible to ask why. One possibility was suggested by Serre \cite{serre2} in his summary of Weil's
mathematical contributions, where he calls attention to the lack of the geometric technology
requisite to a full construction of $Bun_n$, which was subsequently developed only in the
60's by Mumford, Narasimhan, Seshadri, and others \cite{mumford, NS}. However, even with geometric
invariant theory and its relation to $\pi_1$ completed in the remarkable
work of Carlos Simpson \cite{simpson}, there has never been any direct applications of these
moduli spaces (or their cotangent bundles) to Diophantine problems. It is for this reason that the author
locates the difficulty in a far more elementary source, namely, {\em the lack of an
Albanese map to go with $Bun_n$}. Unless  $n=1$, there is no canonical relation\footnote{ It is conceivable
that the theory of Hecke correspondences can be employed  to establish the link.}  between
$Bun_n$ and the points on $X$.
It is  fortunate then that the \'etale topology manages to provide  us with
two valuable tools, namely,  topological fundamental groups that come with fields of definition;
and topological classifying spaces  with extremely canonical
Albanese maps. We owe this to a
 distinguished feature of Grothendieck's theory: the flexible use of basepoints, which
are allowed to be any geometric point at all. The idea that
Galois groups of a certain sort should be regarded as fundamental groups
is likely to be very old, as Takagi\cite{iyanaga} refers to
Hilbert's preoccupation with  Riemann surfaces
as inspiration for class field theory. Indeed, it is true that
that the fundamental group of a smooth variety $V$
will be isomorphic to
the Galois group $\Gal(k(V)^{nr}/k(V))$
 of a maximal unramified extension $k(V)^{nr}$ of its
function field $k(V)$. However, this isomorphism will be {\em canonical} only when
the basepoint is taken to be a separable closure of $k(V)$ that
contains $k(V)^{nr}$:
$$b:\Spec( k(V)^s)\ra \Spec(k(V)^{nr}) \ra \Spec(k(V))\ra V.$$
Within the Galois group approach, there is little room
for small basepoints that come through rational points, or a study
of variation. In fact, there seems to be no reasonable way to
fit path spaces at all into the field picture. This
could then be described as the precise ingredient missing in the arithmetic theory of
fundamental groups at the time of Weil's paper. Even after the introduction of
moving basepoints, appreciation of their genuine usefulness
appears to have taken some time to develop. A rather common
response is to  pass quickly to invariants or situations
where the basepoint can be safely ignored. The author for example came to appreciate
the basepoint as a variable only after reading Professor Deligne's paper written in the 80's \cite{deligne} as well as
the papers of Hodge-theorists like Dick Hain \cite{hain}.

One way to visualize path spaces is to consider
a universal (pro-)cover
$$\tilde{\bX} \rTo \bX.$$
The choice of a lifting $\tb\in \tilde{\bX}_b$ turns the pair into
a {\em universal pointed covering space}. The uniqueness then allows
us to descend to $\Q$, while the universal property determines
canonical isomorphisms
$$\tilde{\bX}_x\simeq \Pet,$$
so that the Galois action can be interpreted using the action on fibers\footnote{The difficult
problem of coming to actual grips with this is  that of constructing a cofinal system making
up $\tilde{\bX}$ in a manner that makes the action maximally visible.
Consider $\mathbb{G}_m$ or an elliptic curve.} .
This is one way to see that the map $\k^{na}$ will never send $x\neq b$
to the trivial torsor, that is, a torsor with an element fixed by $G$, since, by the Mordell-Weil
theorem, nothing but the basepoint will
lift rationally even up to the maximal abelian quotient of $\tilde{\bX}$.
A change of basepoint\footnote{One needs here the elementary
fact that an isomorphism of torsors
$$\pi_1(\bX;b,x)\simeq \pi_1(\bX;b,y)$$
is necessarily induced by a path
$F_x\simeq F_y$.} then shows that the map must in fact be
injective. That is, we have arrived at the striking fact that
points can really be distinguished through the
associated torsors\footnote{It is however, quite interesting to work out
injectivity or its failure for quotients of  fundamental groups corresponding
to other natural systems, like modular towers. Alternatively, one could use the
full fundamental group for a variety where the answer is much less obvious,
like a moduli space of curves.}. In elementary topology, one encounters already the
warning that such path spaces are isomorphic, but not in a canonical fashion. The distinction
may appear pedantic  until one meets such enriched situations
as to endow the torsors with the extra structure
necessary to make them genuinely different.

The remarkable {\em section conjecture} of Grothendieck \cite{grothendieck}
proposes that $\k^{na}$ is even surjective:
$$ X(\Q) \simeq H^1(G,\piet),$$
that is,
\bq
{\em every torsor should be a path torsor.}
\eq
The reader is urged to compare this conjecture with the assertion
that the map
$$\widehat{E(\Q)}\simeq H^1_f(G, \pi^{et}_1(\bE, e)),$$
from Kummer theory is supposed to be bijective for an elliptic curve $(E,e)$. A small
difference has to do with the local `Selmer' conditions on cohomology
indicated by the subscript `$f$', which the complexity of the
non-abelian fundamental group is supposed to render unnecessary. This is
a subtle point on which the experts seem not to offer a consensus. Nevertheless,
the comparison should make it clear to the newcomer that a resolution of
the section conjecture is quite unlikely to be straightforward, being, as it is,
a deep non-abelian incarnation of the principle that suitable conditions
on a Galois-theoretic construction should force it to `come from geometry'\footnote{This notion  in abelian settings coincides roughly with `motivic.'}.
And then, the role of this bijection in the descent algorithm for
elliptic curves might suggest a useful Diophantine context for
the section conjecture \cite{kim3}. Yet another reason for thinking the analogy
through is a hope that the few decades worth of
effort that went into the study of Selmer groups of elliptic curves
might illuminate certain aspects of the section conjecture as well, even at the level of
concrete technology.
\medskip

\medskip

\begin{center}
*
\end{center}
\medskip

\medskip

Our main concern today is with a version of these ideas   where the parallel with elliptic
curves is
especially compelling, in that a good deal of unity between the abelian and non-abelian
realms is substantially realized. This is when
the profinite fundamental group is replaced by the motivic one \cite{deligne}: $$\pi_1^M(\bX,b).$$  The motivic fundamental group lies between the profinite $\pi_1$ and
homology in complexity:
$$\begin{array}{c}
\hat{\pi}_1(\bX,b)\\
| \\
\pi_1^{M}(\bX,b)\\
|\\
H_1(\bar{X})
\end{array}
$$
although it should be acknowledged right away that it is much closer to the bottom
of the hierarchy. The precise meaning of `motivic' should not worry us here more than
in other semi-formal expositions on the subject, since
we will regress quickly to the rather precise use of realizations.
But still, some inspiration may be gather by the rather ghostly presence of
a classifying space
$$H^1_M(G, \pi_1^M(\bX,b))$$
of motivic torsors
as well that of a motivic Albanese map
$$\k^M: X(\Q)\rTo H^1_M(G,\pi_1^M(\bX,b))$$
that associates to points motivic torsors $$\pi_1^M(\bX;b,x)$$
of paths. The astute reader will object that we are again using the
points of $X$ to parametrize motives as in the  subterfuge of
Parshin and Frey, to which we reply that the current family is entirely intrinsic
to the curve $X$, and requires no particular ingenuity to consider.

When it comes to precise definitions \cite{kim1, kim2}, that we must (alas!) inflict upon the reader in a rapid succession
of mildly technical paragraphs, the most important (Tannakian)
category
$$\Un(\bX,\Q_p)$$
consists of  locally constant unipotent $\Q_p$-sheaves on $\bX$,
 where a sheaf is unipotent if it can be constructed using successive extensions starting from
the constant sheaf $[\Q_p]_{\bX}$. As in the profinite theory,
we have a fiber functor
$$F_b:\Un(\bX, \Q_p) \ra \Vect_{\Q_p}$$
that associates to a sheaf $\cV$ its stalk
$\cV_b$, which  has now acquired a linear nature. The $\Q_p$-pro-unipotent
\'etale fundamental group is defined to be
$$U:=\pi_1^{u,\Q_p}(\bX,b):=\Aut^{\otimes}(F_b),$$
the tensor-compatible\footnote{To see the significance
of this notion, one should consider the group
algebra $\C[G]$ of a finite group $G$. On the category $Rep_G(\C)$
of $G$-representations on complex vector spaces, we have the fiber functor
that forgets the $G$-action. Any unit in $\C[G]$ defines an automorphism of
this functor, while the elements of $G$ will then be picked out by the
condition of being tensor-compatible.}
 automorphisms of the fiber functor,
which the linearity equips with the added structure of a pro-algebraic pro-unipotent group over $\Q_p$.
In fact, the descending central series filtration
$$U=U^1\supset U^2 \supset U^3\supset \cdots$$
yields the finite-dimensional algebraic quotients
$$U_n=U^{n+1}\backslash U,$$
at the very bottom of which is an identification
$$U_1=H^{et}_1(\bX,\Q_p)=V_pJ:=T_pJ\otimes \Q_p$$
with the $\Q_p$-Tate module of the (abelian) Jacobian $J$ of $X$.
The different levels are connected by exact sequences
$$0\ra U^{n+1}\backslash U^n \ra U_n \ra U_{n-1} \ra 0$$
that add the extra term $U^{n+1}\backslash U^n$
at each stage, which, however, is a vector group that
can be approached with rather conventional techniques.
In fact, the $G$-action on $U$ lifts the well-studied one on  $V=V_pJ$, and
repeated commutators come together to a quotient map
$$V^{\otimes n} \rOnto U^{n+1}\backslash U^n,$$
placing the associated graded pieces into the category of motives generated
by $J$.
The inductive pattern of these exact sequences
is  instrumental   in making  the unipotent completions considerably more
tractable than their profinite ancestors.

We will again denote by $H^1(G, U_n)$ continuous Galois cohomology with values in the points of
$U_n$. For $n\geq 2$, this is still non-abelian cohomology, and hence, lacks the
structure of a group. Nevertheless, the proximity to homology is evidenced
in the presence of
a remarkable subspace
 $$H^1_f(G,U_n)\subset H^1(G,U_n)$$  defined by local `Selmer' conditions\footnote{Starting
at this point, one should take $p$ to be
a prime of good reduction for $X$,  even though an extension of the theory to the general case should be
straightforward.}
that require the classes to be
\bq
(a)  unramified outside $T=S\cup \{p\}$, where
$S$ is the set of  primes of bad reduction;

(b) and  {\em crystalline} at $p$, a condition coming from $p$-adic Hodge theory.
\eq
The locality of the  conditions refers to their focus on the pull-back of a torsor  for $U$
to the completed fields $\Spec(\Q_l)$. For $l\notin T$, (a) requires the torsor
to trivialize over an unramified extension of $\Q_l$, while
condition (b) requires it to trivialize over Fontaine's ring $B_{cr}$
of crystalline periods \cite{fontaine}. One could equivalently describe the relevant  torsors as
having coordinate rings that are unramified or crystalline as representations
of the local Galois groups.

 Quite important to our purposes is the {\em algebraicity} of the system
$$\cdots \ra H^1_f(G,U_{n+1}) \ra H^1_f(G,U_{n}) \ra H^1_f(G,U_{n-1})\ra \cdots.$$
This is  the {\em Selmer variety} of $X$.
That is, each $H^1_f(G,U_{n})$ is an algebraic variety over $\Q_p$ and the transition maps
are algebraic, so that
 $$H^1_f(G,U)=\{ H^1_f(G,U_n)\}$$ is now a moduli space very similar to  the ones
 that come up in the study of Riemann surfaces \cite{GM}, in that it parametrizes
 crystalline principal bundles for
$U$ in the \'etale topology of $\Spec(\Z[1/S])$. By comparison $H^1(G,\piet)$
has no apparent structure but that of a pro-finite space: the motivic context has
restored some geometry\footnote{`Coefficient geometry,' one might say,
in contrast to $Bun_n$, which carries the algebraic geometry of
the field of definition.} to the  moduli spaces of interest.
The algebraic structure is best understood in terms of $G_T=\Gal(\Q_T/\Q)$,
where $\Q_T$  is the maximal extension of $\Q$ unramified outside $T$. Our moduli space
 $H^1_f(G,U_n)$ sits inside $H^1(G_T, U_n)$ as a subvariety defined by the additional
 crystalline condition.
For the latter, there are sequences
$$0\ra H^1(G_T,U^{n+1}\backslash U^n)\ra H^1(G_T,U_n) \ra H^1(G_T,U_{n-1})\stackrel{\d_{n-1}}{\ra}$$
 $$H^2(G_T,U^{n+1}\backslash U^n)$$
exact in a natural sense, and
the algebraic structures are built up iteratively from the
$\Q_p$-linear structure on the
$$H^i(G_T, U^{n+1}\backslash U^n)$$ using the fact that the boundary maps $\d_{n-1}$ are algebraic\footnote{The reader is warned that it is non-linear in general.}.
That is, $H^1(G_T,U_n)$ is inductively realized as a torsor for the vector group
$H^1(G_T, U^{n+1}\backslash U^n)$
lying over the kernel of $\d_{n-1}$.

It should comes as no surprise at this point that there is a map
$$\bd \k^{u}=\{\k^u_n\}: X(\Q)& \rTo & H^1_f(G,U)\ed $$
associating to a point $x$ the principal $U$-bundle
$$P(x)=\Pux:=\Isom^{\otimes}(F_b,F_x)$$
of tensor-compatible isomorphisms from $F_b$ to $F_x$, that is,
the $\Q_p$-pro-unipotent \'etale paths  from $b$ to $x$. This map is best viewed as
a tower:
$$
\begin{diagram}[height=1.5em]
 &  & \vdots &\\
  & \vdots &  &H^1_f(G,U_4)& \\
  &\ruTo^{\k^u_4}(2,6)&  &\dTo& \\
 & & &H^1_f(G,U_3)&\\
 &\ruTo^{\k^u_3}(2,4) & &\dTo&\\
 &  & &H^1_f(G,U_2) &\\
 & \ruTo^{\k^u_2}  & &\dTo&\\
X(\Q) & \rTo^{\k^u_1}& &H^1_f(G,U_1)&=H^1_f(G,T_p  \otimes \Q_p). \\
\end{diagram}$$
For $n=1$,
$$\k^u_1:X(\Q) \ra H^1_f(G,U_1)=H^1_f(G,T_pJ\otimes \Q_p)$$
reduces to the map from Kummer theory. But the maps
$\k^u_n$ for $n\geq 2$,  much weaker as they are than the $\k^{na}$ discussed
in the profinite context, still do not extend to cycles in any natural way,
and hence, retain   the possibility of separating the structure\footnote{It might be suggested, only half in jest,
that the Jacobian, introduced  by Weil to aid
in the Diophantine study of a curve, has been getting in the way ever since.}  of
$X(\Q)$ from that of $J_X(\Q)$.

 Restricting $U$ to the \'etale site of $\Q_p$, there are local analogues
$$\k^{u}_{p,n}:X(\Q_p) \ra H^1_f(G_p,U_n)$$
that can be  described explicitly (and rather surprisingly) using non-abelian $p$-adic Hodge
theory. More precisely, there is a compatible family of isomorphisms
$$D:H^1_f(G_p, U_n) \simeq U^{DR}_n/F^0$$
to  homogeneous spaces for  the
{\em De Rham fundamental group} $$U^{DR}=\pi_1^{DR}(X\otimes \Q_p,b)$$ of $X\otimes \Q_p$.
Here, $U^{DR}$ classifies unipotent vector bundles with flat connections on $X\otimes \Q_p$,
while $$U^{DR}/F^0$$ is a moduli space for  $U^{DR}$-torsor that carry compatible Hodge filtrations
and Frobenius actions, the latter being obtained from a
 comparison isomorphism\footnote{That is to say,
 if $\cX$ denotes a smooth and proper $\Z_p$-model of $X\otimes \Q_p$,
 the category of unipotent vector bundles with
 flat connections on $X\otimes \Q_p$ is equivalent
 to the category of unipotent convergent isocystals on $\cX\otimes_{\Z_p}\F_p$. This comparison
 is the crucial ingredient in defining $p$-adic iterated integrals \cite{furusho}.} with the crystalline fundamental group and path torsors
 associated to a reduction modulo $p$. The advantage of the De Rham realization is its expression as a $p$-adic
homogenous space whose form is far more transparent than that of Galois cohomology.
The map $D$ (for Dieudonn\'e, as in the theory of $p$-divisible groups)
associates to a crystalline principal bundle $P=\Spec(\cP)$ for $U$, the space
$$D(P)=\Spec([\cP\otimes B_{cr}]^{G_p}).$$ This ends up as a
$U^{DR}$-torsor with Frobenius action and Hodge filtration inherited from that
of $B_{cr}$.
The compatibility of the two constructions is expressed by a diagram
$$\bd
X(\Q_p)& \rTo^{\k^{na}_p}& H^1_f(G_p,U) \\
 & \rdTo^{\k^{u}_{dr/cr}}& \dTo^D\\
 & & U^{DR}/F^0
 \ed$$
 whose commutativity amounts to the non-abelian comparison isomorphism \cite{olsson}
 $$\pi_1^{DR}(X\otimes\Q_p;b,x)\otimes B_{cr}\simeq \Pux \otimes B_{cr}.$$
The explicit nature of the map
 $$\k^{u}_{dr/cr}:X(\Q_p) \ra U^{DR}/F^0,$$
 is a consequence of the $p$-adic iterated integrals\footnote{Special values of
 such integrals have attracted attention because of the connection to values of  $L$-functions.
 Here we are interested primarily in the integrals themselves as analytic functions,
 and in their zeros.} \cite{furusho}
 $$\int_b^z \a_1 \a_2 \cdots \a_n$$
 that appear in its coordinates.
 This expression endows the map  with a highly transcendental nature: for any residue disk $]y[\subset X(\Q_p)$,
 $$\k^{u}_{dr/cr,n}(]y[)\subset U^{DR}_n/F^0$$
 is Zariski dense for each $n$, and is made up  of non-zero convergent power
 series that are obtained explicitly as repeated anti-derivatives starting from differential forms on $X$.

 Finally, the local and global constructions fit into a family of commutative diagrams
 $$\bd
 X(\Q) & \rTo& X(\Q_p)& &\\
 \dTo& & \dTo& \rdTo& \\
 H^1_f(G,U_n)& \rTo^{\loc_p}&H^1_f(G_p,U_n)& \rTo^D& U^{DR}_n/F^0
 \ed$$
 where the bottom horizontal maps are algebraic and the vertical maps
 transcendental. Thus, the difficult inclusion
 $X(\Q)\subset X(\Q_p)$ has been replaced by the map\footnote{The strange notation
 is comes the view that $D$ is itself a log map, according to Bloch and Kato \cite{BK}.}
 $\log_p:=D\circ \loc_p$, whose algebraicity gives a glimmer of hope that the arithmetic geometry
 can be understood and controlled.

The following result is basic to the theory.
\begin{thm}
Suppose
$$ \log_p(H^1_f(G,U_n))\subset U^{DR}_n/F^0$$
is not Zariski dense
for some $n$. Then
$X(\Q)$  is finite\footnote{Professor Serre would object that the theorem is trivially
true since $X(\Q)$ {\em is} finite. The author offers no defense.}.
\end{thm}

The proof of this assertion in its entirety is captured by the diagram
$$\begin{diagram}X(\Q) & \rInto &X(\Q_p)  \\
\dTo^{\k^{u}_n} & & \dTo^{\k^{u}_{dr/cr,n} }\\
H^1_f(G, U_n) &  \rTo^{\log_p} & U^{DR}_n/F^0 \\
 & &   \dTo_{\exists \phi \neq 0 } \\
 & &   \Q_p
\end{diagram}$$
indicating the existence of a non-zero algebraic function
$\phi$ vanishing on $\log_p(H^1_f(G, U_n))$. Hence, the function
$\phi \circ \k^{u}_{dr/cr,n}$ on $X(\Q_p)$ vanishes on
$X(\Q)$. But
this function is a non-vanishing
convergent power series on each residue disk, which therefore can have only finitely many zeros.
$\Box$.

A slightly more geometric account of the proof might point to the fact that the image of $X(\Q_p)$
in $U^{DR}_n/F^0$
is a space-filling curve, with no portion contained in a proper subspace.
Hence, its intersection with any proper subvariety must be discrete. Being compact as
well, it must then be finite\footnote{This proof, involving a straightforward interplay
of denseness, non-denseness, and compactness, is a curious avatar of
some ideas of Professor Deligne relating the section conjecture to Diophantine finiteness. }.
Serge Lang once proposed a strategy for proving the Mordell conjecture by deducing it from
a purely geometric hope that the complex points on a curve of higher genus
might intersect a finitely generated subgroup of the Jacobian in finitely many points.
While that idea turned out to be very difficult to realize, here we have a non-Archimedean
analogue, wtih $U^{DR}_n/F^0$ playing the role of the complex  Jacobian,
and the Selmer variety that of the Mordell-Weil group.

The hypothesis of the theorem on non-denseness of the global Selmer variety
is expected always to hold for $n$ large, in that we should have \cite{kim2}
$$\dim H^1_f(G,U_n)<< \dim U^{DR}_n/F^0.$$
(Recall that the map $\log_p$ is algebraic.)
Such an inequality follows, for example, from the reasonable folklore conjecture that$$H^1_f(G, M)=0$$
for a motivic Galois representation\footnote{It suffices here to
take $M$ to be among the motives generated by $H^1(X)$.}  $M$ of weight $>0$. This, in turn,
might be deduced from the conjecture of Fontaine and Mazur on Galois representations
of geometric origin \cite{FM}, or from portions of the Bloch-Kato conjecture\footnote{We thus have reason, in the manner of
physicists, to
regard Theorem 1 as good news for mixed motives, in that highly non-trivial real phenomena are among the corollaries of their theory. A small counterpoint
to the pessimistic view of Professor Serre.} \cite{BK}. The point is that
if we recognized the elements of $H^1_f(G, M)$ themselves to be motivic,
then the vanishing would follow from the existence of a weight filtration.
Thus instead of the implication
\bq
Non-abelian `finiteness of $\Sha$' (= {\em section conjecture}) $\Rightarrow$ finiteness of $X(\Q)$.
\eq
expected by Grothendieck, we have
\bq
`Higher abelian finiteness of $\Sha$' (that $H^1_f(G,M)$ is generated by
motives)$\Rightarrow$ finiteness of $X(\Q)$.
\eq
This is not the only place that our considerations revolve around pale shadows of
the section conjecture. One notes for example, the critical use of the
dense image of $\k^{u}_{dr/cr}$, which could itself  be thought of as an `approximate local section
conjecture.'

In spite of all such lucubrations (that fascinate the author and quite likely no one else), we must now face the
 plain and painful fact that an unconditional proof of the hypothesis  for
large $n$ (and hence, a new proof of finiteness) can be given only in situations where the image of $G$ inside
$\Aut(H_1(\bX, \Z_p))$ is {\em essentially abelian}. That is, when
\bq
-$X$ is an affine hyperbolic of genus zero (say $\P^1\setminus\{0,1,\infty\}$) \cite{kim1};

-$X=E\setminus\{e\}$ for an elliptic curve $E$ with complex multiplication \cite{kim4};

-(with John Coates)  $X$ is compact of genus $\geq 2$ and the Jacobian $J$ factors into abelian varieties with potential
complex multiplication \cite{CK}.
\eq
The first two cases require a rather obvious modification tailored to the study of
integral points, while  the two CM cases require  $p$ to be split inside the CM fields.
Given the intermediate state of the purported application,
the reason  for persevering in an abstruse investigation of known results
might seem obscure indeed. We will return to this point
towards the end of the lecture, side-stepping the issue  for now in favor of a brief sketch of the methodology, confining our
attention to the third class of curves.

There is a pleasant quotient\footnote{For $\P^1\setminus \{0,1,\infty\}$,
such quotients came up in the process of isolating (simple-)polylogarithms \cite{BD}.}$$U \rOnto W:=U/[[U,U],[U,U]]$$
of $U$ that allows us to extend the key diagrams.
$$\begin{diagram}X(\Z_S) & \rInto &X(\Z_p)& && &   \\
\dTo^{\k^{u}_n} & & \dTo^{\k^{u}_{p,n} } & \rdTo^{\k^u_{dr/cr,n}}&\\
H^1_f(G, U_n) &  \rTo^{\mbox{loc}_p} & H^1_f(G_p, U_n)& \rTo^{D}&  U^{DR}_n/F^0 \\
\dTo &  & \dTo & & \dTo\\
H^1_f(G, W_n) &  \rTo^{\mbox{loc}_p} & H^1_f(G_p, W_n)&\rTo^{D} & W^{DR}_n/F^0
\end{diagram}$$
The structure of $W$ turns out to be much simpler than that of $U$,
and we obtain the following result.
\begin{thm}[with John Coates]
Suppose $J$ is isogenous to a product of abelian varieties having potential
complex multiplication. Choose the prime $p$ to split in all the $CM$ fields
that occur. Then
$$\dim H^1_f(G,W_n) < \dim W^{DR}_n/F^0$$
for $n$ sufficiently large.
\end{thm}
The non-denseness of $\log_p(H^1_f(G,U))$ is an obvious corollary.

We give an outline of the proof assuming $J$ is simple.
Since $$\dim H^1_f(G,W_n) \leq \dim H^1(G_T,W_n),$$
it suffices to estimate the dimension of cohomology with restricted ramification.
Via the exact sequences
$$0\ra H^1(G_T,W^{n+1}\backslash W^n)\ra H^1(G_T,W_n) \ra H^1(G_T,W_{n-1})$$
the estimate can be reduced to a sum of abelian ones:
$$\dim H^1(G_T,W_n)\leq \sum_{i=1}^n\dim H^1(G_T,W^{i+1}\backslash W^i).$$
The linear representations
$W^{i+1}\backslash W^i$ come with Euler characteristic formulas\footnote{The minus sign in
the superscript refers to the negative eigenspace of complex conjugation.  This has roughly
half the dimension of the total space, and ends up unduly important to our estimates.} \cite{milne}:
$$\dim H^0(G_T,W^{i+1}\backslash W^i)-\dim H^1(G_T,W^{i+1}\backslash W^i)$$
$$+\dim H^2(G_T,W^{i+1}\backslash W^i)=
-\dim [W^{i+1}\backslash W^i]^{-}.$$
out of which the $H^0$ term always vanishes, leaving
$$ \dim H^1(G_T,W^{i+1}\backslash W^i)=\dim [W^{i+1}\backslash W^i]^{-}+ \dim H^2(G_T,W^{i+1}\backslash W^i).$$
The comparison with the topological fundamental group of $X(\C)$ reveals $U$ to
be the unipotent completion of a free group on $2g$ generators modulo a single
relation. This fact can applied to construct a  Hall basis for the Lie algebra of $W$ \cite{reutenauer}, from
which we get an elementary estimate
$$\sum_{i=1}^n \dim [W^{i+1}\backslash W^i]^{-}\leq [(2g-1)/2]\frac{n^{2g}}{(2g)!}+O(n^{2g-1}).$$
Similarly, on the De Rham side the dimension
$$\dim W^{DR}_n/F^0=W_2/F^0+\sum_{i=3}^n\dim[W^{DR,i+1}\backslash W^{DR,i}]$$
can easily be bounded below by
$$(2g-2)\frac{n^{2g}}{(2g)!}+O(n^{2g-1}).$$
Hence, since $g \geq 2$, we have
$$\sum_{i=1}^n \dim [W^{i+1}\backslash W^i]^{-}<< \dim W^{DR}_n/F^0.$$
Therefore, it remains to show that
$$\sum_{i=1}^n \dim H^2(G_T,W^{i+1}\backslash W^i)=O(n^{2g-1}).$$
Standard arguments with Poitou-Tate duality\footnote{which switches the focus
from $H^2$ to $H^1$ at the cost of dealing with some insignificant local terms} \cite{milne} eventually reduces the problem to
the study of $$\Hom_{\G}[M(-1), \sum_{i=1}^n[W^{i+1}\backslash W^i]^*],$$
where
\bq
-$F$ contains  $\Q(J[p])$ and is a field of
definition for all the complex multiplication;

-$\G=\Gal(F_{\infty}/F)$ for the field $$F_{\infty}=F(J[p^{\infty}])$$ generated by
the $p$-power torsion of $J$;

-and $$M=\Gal(H/F_{\infty})$$ is the Galois group of the $p$-Hilbert class field $H$
of $F_{\infty}$.
\eq
Choosing an annihilator\footnote{provided by a theorem of Greenberg \cite{greenberg1}} $$\cL\in \Lambda:=\Z_p[[\G]]\simeq \Z_p[[T_1, T_2, \ldots, T_{2g}]]$$
for $M(-1)$ in the Iwasawa algebra,  we need to count its zeros among the characters
that appear in $$\sum_{i=1}^n[W^{i+1}\backslash W^i]^*.$$
If$\{\psi_i\}_{i=1}^{2g}$ are the characters that make up $H^1(\bX,\Q_p)$,
the characters  in
$[W^{i+1}\backslash W^i]^*$ are a subset of
$$\psi_{j_1}\psi_{j_2}\psi_{j_3}\cdots \psi_{j_{i}},$$
where
$j_1<j_2\geq j_3\geq  \cdots \geq j_i$.
 After a change of variables, a lemma of Greenberg \cite{greenberg2} allows us to assume a form
 $$\cL=a_0(T_1, \ldots, T_{2g-1})+a_1(T_1, \ldots, T_{2g-1})T_{2g}+\cdots$$
$$+
a_{l-1}(T_1, \ldots, T_{2g-1})T_{2g}^{l-1}+T_{2g}^l,$$
a polynomial in $T_{2g}$. We can estimate the number of zeros
by considering instead the $2g-1$ polynomials obtained by fixing
the index $j_1$, and counting their zeros among the set of
$\psi_{j_2}\psi_{j_3}\cdots \psi_{j_{i}}$
with $j_2\geq j_3\geq \cdots \geq j_i$.
As $i$ runs from 1 to $n$, the multi-indices in the exponents of
 $$\psi_1^{m_1}\psi_2^{m_2}\cdots \psi_{2g}^{m_{2g}}$$
 that occur  are among the integer points in a simplex of side length $n-1$ in
a space of dimension $2g$.  But then, since the coefficients $a_i$ depend only
on the projection of these integer points to a simplex of one smaller dimension, and
the number of zeros lying above each such point is at most $l$, the total number
of zeros is $O(n^{2g-1})$.
Since $M(-1)$ is $\L-$finitely-generated, we deduce the
bound
$$\Hom_{\G}[M(-1), \sum_{i=1}^n[W^{i+1}\backslash W^i]^*]=O(n^{2g-1})$$
desired. $\Box$
\medskip

We  have now set up the first genuine occasion to motivate our constructions. The annihilator
$\cL$
is a version of an {\em algebraic $p$-adic $L$-function} controlling the situation.
It is therefore of non-trivial interest that the sparseness of its zeros is responsible for the finiteness
of points.
The parallel with the case of elliptic curves \cite{CW, kato, kolyvagin, rubin} might be seen clearly by
comparing the implications
\bq
non-vanishing of $L$ $\Rightarrow$ control of Selmer groups $\Rightarrow$
finiteness of points
\eq familiar from the arithmetic of elliptic curves to
the one given:
\bq
sparseness of $L$-zeros $\Rightarrow$ control of Selmer {\em varieties} $\Rightarrow$
finiteness of points.
\eq
As promised, the motivic fundamental group has provided a natural thread
linking  abelian and non-abelian Diophantine problems.

We remark that the non-CM case could  proceed along the same
lines, except that the group $\G$ and hence,
the corresponding Iwasawa algebra is non-abelian. But the fact remains
that the estimate
$$\dim \Hom_{\L}(M, \oplus_{i=1}^{n} W^{i+1}\backslash W^i)=O(n^{2g-1})$$
is sufficient for the analogue of Theorem 2, and hence, for the finiteness
of points. The representation $W^{i+1}\backslash W^i$ is a subquotient of
the more familiar one
$$(\Lambda^2 V_p)\otimes (\Sym^{i-2}V_p)$$
and the difference in dimensions is likely to
count for very little in the  coarse estimates. It might therefore be easier to
work with
$$\Hom_{\L}[M, \oplus_{i=1}^{n-2} (\Lambda^2 V_p)\otimes (\Sym^{i}V_p)].$$
Otmar Venjakob \cite{venjakob} has shown that $M$ is locally torsion, so that
a generating set $\{m_1, m_2,\ldots, m_d\}$ for $M$
determines for each $i$ a non-commutative power series $ f_i\in \L$
annihilating  $m_i$.
We  must then count the {\em non-abelian zeros}\footnote{As noted by Mahesh Kakde, it would be nice
to know enough to formulate this in terms of a characteristic element
$f\in K_1(\L_{S^*})$ for $M$, whereby the count will
be of irreducible representations $\rho: \G \ra N$ for which
$f(\rho)=0$ \cite{CFKSV}.}
of $f_i$, that is, the representations containing vectors annihilated by
$f_i$ among the irreducible factors of
$\oplus_{i=1}^{n-2}(\Lambda^2 V_p)\otimes (\Sym^{i}V_p)$.

John Coates has stressed the role played by the ideal class group $M$ in this picture, which
is a priori smaller than the Iwasawa module relevant to elliptic curves.
The reason that ramification at $p$ can be ignored for now is that the local contribution at $p$
is also of lower order as a function of  $n$. For the Diophantine geometry of
abelian fundamental groups, however, the option of passing to large $n$ is absent. One is tempted
 to offer this as a kind of explanation for the infinitely many rational points that can live on an elliptic
curve.
\medskip

\medskip
\begin{center}
*
\end{center}
\medskip

\medskip

Some preliminary evidence at present suggests another reason to pursue
a $\pi_1$ approach to finiteness \cite{kim5}.
This is the possibility that the function $\phi$ occurring in the proof of theorem 1
can be made explicit, leading to analytic defining equations for
$$X(\Q)\subset X(\Q_p).$$
 For one thing, the map
$$\log_p:H^1_f(G,U_n)\ra U^{DR}_n/F^0$$
occurs in the category of algebraic   varieties over $\Q_p$,  and is therefore  amenable (in principle) to
computation \cite{coleman, poonen}. Whenever the map itself can be presented, the computation
of the image is then a matter of applying standard algorithms.
A genuinely {\em feasible} approach, however, should be effected by the {\em cohomological construction} of a function $\psi$ as below that
vanishes on global classes.
$$\bd
 X(\Q) & \rTo& X(\Q_p)& &\\
 \dTo& & \dTo& \rdTo& \\
 H^1_f(G,U_n)& \rTo^{\loc_p}&H^1_f(G_p,U_n)& \rTo^D& U^{DR}_n/F^0\\
 & &\dTo^{\psi} &\ldTo^{\phi} & \\
 & & \Q_p& &
 \ed$$
 That is, once we have  $\psi$, we can put
 $$\phi=\psi\circ D^{-1},$$
a function  whose precise computation might be regarded as a `non-abelian explicit reciprocity law.'
The vanishing itself should be explained by a local-to-global reciprocity, as in the work of Kolyvagin, Rubin, and Kato on the
conjecture of Birch and Swinnerton-Dyer \cite{kolyvagin, rubin, kato}.

These speculations are best given substance with an example, albeit in an affine setting.
Let $X=E\setminus \{e\}$, where
$E$ is an elliptic curve of rank 1 with $\Sha(E)[p^{\infty}]=0$.
The significance of the hypotheses is that the $\Q_p$ localization
map is bijective on points,
$$loc_p:E(\Q)\otimes \Q_p\simeq H^1_f(G_p,V_p(E)),$$
and the second cohomology with restricted ramification vanishes:
$$H^2(G_T, V_p(E))=0.$$
We will construct a diagram:
$$\bd
 X(\Z) & \rTo& X(\Z_p)& &\\
 \dTo& & \dTo& \rdTo& \\
 H^1_{f,\Z}(G,U_2)& \rTo^{\loc_p}&H^1_f(G_p,U_2)& \rTo^D& U^{DR}_2/F^0\\
 & &\dTo^{\psi} &\ldTo^{\phi} & \\
 & & \Q_p.& &
 \ed$$
 using just the first non-abelian level $U_2$
 of the unipotent fundamental group.
 We have introduced here a refined Selmer variety $H^1_{f,\Z}(G,U_2)$ consisting of classes that are
 actually trivial
 at all places $l\neq p$. It is a relatively straightforward matter to show  that the integral points land in this
 subspace \cite{KT}.

The relevant structure now is  a Heisenberg group
$$0\ra \Q_p(1)\ra U_2 \ra V_p \ra 0,$$
that we will analyze in terms of the corresponding extension of Lie algebras
$$0\ra \Q_p(1)\ra L_2 \ra V_p \ra 0.$$
Conveniently, at this level, the Galois action on $L_2$ splits\footnote{This uses
the multiplication by $[-1]$, as in Mumford's theory of theta functions.}:
$$L_2=V_p\oplus \Q_p(1),$$
provided we use  a tangential base-point at the missing point $e$.
With the identification\footnote{For unipotent groups, the power series
for the log map stops after finitely many terms,
defining an algebraic isomorphism. The group then can be
thought of as the Lie algebra itself with a twisted binary operation given by the
Baker-Campbell-Hausdorff formula \cite{serrelie}.} of $U_2$ and $L_2$,
non-abelian cochains can be thought of as maps
$$\xi:G_p\rTo L_2$$
and expressed in terms of components $\xi=(\xi_1,\xi_2)$
with respect to the decomposition. The cocycle condition in these coordinates reads\footnote{In his book on
gerbes, Breen emphasizes the importance of a familiarity with the `calculus of cochains.' Indeed, the typical
number-theorist will be quite anxious about non-closed cochains like $\xi_2$. Unfortunately,
they are as unavoidable as the components of connection forms in non-abelian gauge theory, which
obey complicated equations even when
the connections themselves are closed in a suitable sense.}
$$d\xi_1=0, \ \ \ \ \ d\xi_2=(-1/2)[\xi_1,\xi_1].$$
Define
$$\psi(\xi):=[\loc_p(x),\xi_1]-2\log \chi_p\cup \xi_2\in H^2(G_p, \Q_p(1))\simeq \Q_p,$$
where
$$\log \chi_p:G_p\ra \Q_p$$ is the logarithm of the $\Q_p$-cyclotomic character and
$x$ is a {\em global } solution,
to the equation
$$dx=\log \chi_p \cup \xi_1.$$
The equation makes sense on $G_T$ since both $\chi_p$ and $\xi_1$ have
natural extensions to global classes, while the non-trivial   existence of
the global solution $$x:G_T\ra V_p$$ is guaranteed by the aforementioned  vanishing  of $H^2$.
One checks readily that $\psi(\xi)$ is indeed a 2-cocycle whose class is independent of the choice
of $x$.
\begin{thm}
$\psi$ vanishes on the image
of
$$\loc_p: H^1_{f,\Z}(G,U_2)\ra H^1_f(G_p, U_2).$$
\end{thm}
The proof is a simple consequence of the standard reciprocity sequence
$$0\ra H^2(G_T,\Q_p(1))\ra \oplus_{v\in T} H^2(G_v,\Q_p(1))\ra \Q_p\ra 0.$$
The point is that if $\xi$ is global then so is $\psi(\xi)$. But this class has
been constructed to vanish at all places $l\neq p$. Hence, it must also vanish at $p$.

An explicit formula on the
De Rham side in this case is rather  easily obtained. Choose a Weierstrass equation for $E$ and let
$$\a=dx/y, \ \ \  \b=xdx/y.$$
Define
$$\log_{\a}(z):=\int_b^z\a, \ \ \ \log_{\b}(z):=\int_b^z \b,$$
$$D_2(z):=\int_b^z \a \b,$$
via (iterated) Coleman integration.
\begin{cor}
For any two points  $y,z\in X(\Z)\subset X(\Z_p)$, we have
$$\log^2_{\a}(y)(D_2(z)-\log_{\a}(z)\log_{\b}(z))=\log^2_{\a}(z)(D_2(y)-\log_{\a}(y)\log_{\b}(y)).$$

\end{cor}
The proof  uses an action of the multiplicative monoid
$\Q_p$ on $H^1_f(G,U_2)$ that covers the scalar multiplication
on $E(\Q)\otimes \Q_p$. That is,
$$\lambda \cdot (\xi_1, \xi_2)=(\lambda \xi_1, \lambda^2 \xi_2).$$
Evaluating $\psi$ on the class
$$\log_{\a}(x) \k^u_2(y)-\log_{\a}(y)\k^u_2(x)\in H^1_f(G_p, U^3\backslash U^2)$$
leads directly to the formula displayed. The harmonious form of the resulting
constraint is perhaps an excuse for some general optimism.
Of course, as it stands, the formula is useful only if there is a point $y$ of infinite order already at hand.
One can then look for the other integrals points in the zero set of
the function
$$D_2(z)-\log_{\a}(z)\log_{\b}(z)-(\frac{D_2(y)-\log_{\a}(y)\log_{\b}(y)}{\log^2_{\a}(y)})\log^2_{\a}(z)$$
in the coordinate $z$.

The {\em meaning} of the construction given is not yet clear to the author,
even as some tentative avenues of interpretation are opening  up quite recently. If the analogy with
the abelian case is to be taken seriously, $\psi$ should be a small fragment
of {\em non-abelian duality} in Galois cohomology\footnote{Kazuya Kato's immediate reaction
to the idea of  non-abelian duality was that
it should have an `automorphic' nature.  Such a suggestion might be highly relevant if the {\em reductive
completion} of fundamental groups could somehow be employed in an arithmetic setting.
For the unipotent completions under discussion, the author's inclination is to look
for duality that is a relatively straightforward lift of the abelian phenomenon.}.
For the abelian quotient, one has the usual duality
$$H^1(G_p,V) \times H^1(G_p, V^*(1)) \rTo H^2(G_p, \Q_p(1))\simeq \Q_p$$
with respect to which $H^1_f(G_p,V)$ and $H^1_f(G_p, V^*(1))$
are mutual annihilators. We take the view that $$H^1(G_p,V^*(1))/H^1_f(G_p,V^*(1))$$
is thereby  a systematic source of functions on $H^1_f(G_p, V)$,
which  can  then be used to annihilate global classes when the function itself
comes from a suitable class\footnote{The author is not competent to review
here the laborious procedure for producing such classes as was developed
in the work of Kolyvagin and Kato. The guiding concept in the abelian case is that
of a {\em zeta element}.} in $H^1(G, V^*(1))$. After a minimal
amount of non-commutativity has been introduced, our $\psi$ is exactly such
a global function on the local cohomology $H^1_f(G_p,U_2)$ that ends up thereby {\em annihilating the Selmer variety}.
The main difficulty is that we know not yet a suitable  space in which $\psi$ lives.
Allowing ourselves a further flight of fancy, the elusive function in general
might eventually be the subject of an Iwawasa theory rising out of a landscape
radically more non-abelian and non-linear than we have dared to dream of thus far \cite{kato2}.

\medskip

\medskip
\begin{center}
*
\end{center}
\medskip

\medskip
It has been remarked that the title of this lecture was chosen to be maximally
ambiguous. Notice, however, that Galois theory in dimension zero, according to Galois, proposes
groups as  structures encoding the
Diophantine geometry of equations in one variable.
The proper subject of Galois theory in dimension one  should then be
a unified network of structures relevant to
the Diophantine geometry of polynomials in two variables.
Included therein one may find the arithmetic fundamental groups,
motivic $L$-functions of weight one, and moduli spaces of torsors
that have already proved their scattered usefulness to the trade\footnote{The section conjecture says the set of points on
a curve of higher genus {\em is}
a moduli space of torsors. One might take this to
be a categorical structure that generalizes the abelian groups
that come up in elliptic curves.}. The picture as
a whole is blatantly far from clear, coherent, or complete at this stage
\footnote{It has been an enduring  source of
amazement to the author that true number-theorists
employ philosophies that never work in practice as planned at the outset.
The numerous subtle twists and turns that one may find, for example,
in the beautiful theorems of Richard Taylor, that adhere nevertheless to
the overall  form
of a grand plan, are hallmarks of the
kind of artistry that a mere generalist could never aspire to. It is essentially
for this reason that the author has avoided thus far the question of applying the techniques of
this paper to varieties of
higher dimension, for example, those with a strong degree of hyperbolicity.
A theory whose end product is a single function
applies immediately only in dimension one. It is not inconceivable
that an arsenal of clever tricks
will strengthen the machinery shown here to make
it more broadly serviceable. A robust
strategy that makes minimal demands on the user's ingenuity, however,
should  expect the requisite structures to evolve as one
climbs up the dimension ladder, perhaps in a manner reminiscent of
Grothendieck's {\em poursuite}.}.

{\footnotesize  Department of Mathematics, University College London,
Gower Street, London, WC1E 6BT, United Kingdom and The Korea Institute
for Advanced Study, Hoegiro 87, Dongdaemun-gu, Seoul 130-722, Korea}

\end{document}